\documentclass[12pt]{article}
\usepackage[centertags]{amsmath}

\newtheorem{proposition}{Proposition}
\newtheorem{lemma}[proposition]{Lemma}
\newtheorem{remark}[proposition]{Remark}

\def\A{{\cal A}}

\def\bcn{q\,\nu\,b+c}
\def\C{{\cal C}}
\def\Co{{\bf C}}
\def\De{\Delta}
\def\fidi{\hskip5pt \vrule height4pt width4pt depth0pt \par}
\def\fraz#1#2{{\strut\displaystyle #1\over\displaystyle #2}}
\def\half{{\scriptstyle{\frac1 2}}}
\def\K{{\cal K}}
\def\L{{\cal L}}
\def\Na{{\bf N}}
\def\R{{\cal R}}
\def\Re{{\bf R}}
\def\slr{SL(2,\Re)}
\def\slqr{SL_q(2,\Re)}
\def\tens{\otimes}

\def\vep{\varepsilon}
\def\Ze{{\bf Z}}

\begin{document}

\title{The Coisotropic Subgroup Structure of
$\slqr$}

\author{F.Bonechi$\,{}^{1,2}$, N.Ciccoli$\,{}^{3}$, R.Giachetti$\,{}^{2,1}$,
E.Sorace$\,{}^{1,2}$, M.Tarlini$\,{}^{1,2}$}

\date{June 1999}

\maketitle

\centerline{{\small  ${ }^1$ INFN Sezione di Firenze}}
\centerline{{\small ${ }^2$ Dipartimento di Fisica, Universit\`a di
Firenze, Italy. }}
\centerline{{\small ${ }^3$ Dipartimento di Matematica, Universit\`a
di Perugia, Italy. }}

\begin{abstract}
\noindent
{We study the coisotropic subgroup structure of standard 
$\slqr$ and the corresponding embeddable quantum homogeneous
spaces. While the subgroups ${\bf S}^1$ and $\Re_+$ survive undeformed 
in the quantization as coalgebras, we show that $\Re$ is deformed to a 
family of quantum coisotropic subgroups whose coalgebra can not be 
extended to an Hopf algebra. We explicitly describe the quantum homogeneous
spaces and their double cosets.}
\end{abstract}

\medskip
\noindent
{\bf Math.Subj.Classification}: 16W30, 17B37, 81R50

\smallskip
\noindent
{\bf Keywords}: $\slqr$, coisotropic quantum subgroups, quantum
homogeneous spaces.

\thispagestyle{empty}

\bigskip\bigskip
\noindent
{\bf 1.} The purpose of this paper is the study of the subgroup
structure of standard $\slqr$, together with its embeddable quantum
homogeneous spaces. Let us first recall the classical case. A 
geometrical picture can be given by considering 
the adjoint action of $\slr$ on its Lie algebra $sl(2,\Re)$. 
The Killing form, invariant with respect to the adjoint action, has 
signature $(1,2)$. After the identification of $sl(2,\Re)$ with
$\Re^3$, one obtains that the group action preserves all quadratic 
submanifolds of the form $x^2-y^2-z^2=\Theta\,$, with $\Theta\in\Re$. 
Excluding the trivial orbit ${(0,0,0)}$ and studying the isotropy
subgroups of points of such quadrics, one can distinguish three
essentially different cases of homogeneous spaces:
\begin{itemize}
\item[1.] When $\Theta<0$ one gets one-sheeted hyperboloids; the
corresponding isotropy subgroups are all conjugated to
$A=\left\{(\begin{smallmatrix}\lambda & 0\\ 0 & \lambda^{-1} 
\end{smallmatrix})\;,\lambda\in\Re_+\right\}$.
These subgroups can also be characterized as those containing
hyperbolic matrices, {\it i.e.} matrices $M$ such that $|{\rm tr} M|>2$. We
will call such subgroups of $\Re_+$-type.
\item[2.] When $\Theta>0$ one gets two sheeted hyperboloids. 
The connected components of the isotropy subgroups are conjugated to 
$K=\left\{(\begin{smallmatrix}\cos\theta & 
\sin\theta\\  -\sin\theta & \cos\theta \end{smallmatrix})\;,
\theta\in[0,2\pi)\right\}$. 
These subgroups contain elliptic matrices, {\it i.e.} matrices $M$ 
such that $|{\rm tr} M|<2$. We will call such subgroups of ${\bf S}^1$-type.
\item[3.] Finally when $\Theta=0$ one gets the light-cone minus its
vertex. The connected components of the isotropy subgroups are
conjugated to $N=\left\{(\begin{smallmatrix} 1 & b\\ 0 & 1 
\end{smallmatrix})\;,b\in \Re\right\}$. All of its elements are
parabolic, {\it i.e.} matrices $M$ such that $|{\rm tr} M|=2$.
We will call such subgroups of $\Re$-type.

\end{itemize}
All one-dimensional subgroups of $\slr$ belong to one of the conjugacy 
classes of the isotropy subgroups above. The representatives we have
chosen are those appearing in the Iwasawa decomposition $\slr=K\,A\,N\,.$

Consider now the quantum situation. Quantum subgroups (quotients by
Hopf ideals) do not play the same role as subgroups in the
classical case. Indeed it is well known (see {\it e.g.}  
\cite{P1}) that in the quantum case the limited number of subgroups
that survive is not enough to obtain all the homogeneous spaces. 

Recently, however, the concept of quantum subgroup has been
generalized to that of coisotropic quantum subgroup, allowing to
recover from a quotient type procedure all known families
of quantum embeddable homogeneous spaces (see for example \cite{Br} 
for quantum spheres, and \cite{Ci1} for quantum planes and
cylinders). Coisotropic subgroups are quotient by a coideal, 
right (or left) ideal so that they inherit only the coalgebra 
structure while the algebra is weakened to a right (or left) module. 
Although the involution doesn't pass to the quotient nevertheless the 
real structure survives with the definition of $\tau=*\circ S$
\cite{BCGST}. Indeed, according to this structure, we are able to 
describe the different nature of coisotropic subgroups of the 
standard $\slqr$ with respect to the ones defined by $SU_q(2)$
fibrations over quantum-spheres.  In addition to subgroups of ${\bf S}^1$ 
and $\Re_+$ type that remain undeformed as coalgebras, we find 
a family of coisotropic quantum subgroups whose coalgebras are all
isomorphic to a coalgebra $\Re_q$. The latter deforms the subgroups 
of $\Re$ type and cannot be completed to a Hopf algebra. 
The corresponding homogeneous spaces are the 
analogues of the exceptional quantum-spheres parametrized in \cite{P2} 
by $c=c(n)$. Contrary to these  exceptional quantum-spheres our
special series  have classical points so that they are embeddable. 
In Proposition \ref{Re_q} we give a detailed analysis of $\Re_q$
showing that it is not cosemisimple. 

We construct the embeddable quantum homogeneous spaces as the spaces
of coinvariant elements and we recover the results given in \cite{Lei} 
starting from the semiclassical covariant Poisson structure.
By the use of coisotropic subgroups we are able to define the 
double cosets and give their explicit description.
Moreover the properties of these homogeneous spaces are clarified by 
our geometrical construction.

\bigskip\bigskip
\noindent
{\bf 2.} In this section we recall the definition and the main
properties of real coisotropic quantum subgroups as well as their 
associated homogeneous spaces (see \cite{Ci1,BCGST} for more details). 
We then discuss the equivalence of coisotropic subgroups determined by 
the characters of the whole Hopf algebra. 

Given a real quantum group $(\A,*)$ we will call {\it real
coisotropic quantum right $($left$)$ subgroup} $(\K, \tau_\K)$ a
coalgebra, right (left) $\A$-module $\K$ such that:
\begin{itemize}
\item[i)] there exists a surjective linear map $\pi:
\A\rightarrow \K$, which is a morphism of coalgebras and of 
$\A$-modules (where $\A$ is considered as a module on itself 
via multiplication);
\item[ii)] there exists an antilinear map $\tau_\K:
\K\rightarrow \K$ such that
$\tau_\K\circ\pi=\pi\circ\tau$, where $\tau = *\circ S$.
\end{itemize}
A $*$-Hopf algebra ${\cal S}$ is said to be a real quantum subgroup 
of $\A$ if there exists a $*$-Hopf algebra epimorphism 
$\pi:\A\rightarrow {\cal S}$; evidently this is a particular
coisotropic subgroup.
We remark that a coisotropic quantum subgroup is not in general 
a $*$-coalgebra but it has only $\tau_\K$ defined on it. It is easy 
to verify that if a coisotropic quantum subgroup is also a
$*$-coalgebra and $\pi\circ *=*\circ\pi$, it is then possible to
complete the structure so to have a quantum subgroup.

Coisotropic quantum subgroups are characterized by the
following proposition.
\medskip
\begin{proposition}
\label{pro_coiso}
There exists a bijective correspondence between coisotropic quantum
right $($left$)$ subgroups and $\tau$-invariant two-sided coideals, right
$($left$)$ ideals in $\A$. \fidi
\end{proposition}
\medskip
A canonical construction for embeddable quantum homogeneous spaces can 
now be provided.
\medskip
\begin{proposition}
If $(\K,\tau_\K)$ is a right coisotropic quantum
subgroup of $(\A,*)$ then
$$
B^\pi = \{a\in\A \,|\, ({\rm id}\otimes\pi)\Delta a = a\otimes\pi(1)\}
$$
is a right embeddable quantum homogeneous space.

If $(\K,\tau_\K)$ is a left coisotropic quantum subgroup of
$(\A, *)$ then
$$
B_\pi = \{a\in\A \,|\, (\pi\otimes{\rm id})\Delta a = \pi(1)\otimes a\}
$$
is a left embeddable quantum homogeneous space. \fidi
\end{proposition}
\medskip
The correspondence between coisotropic quantum subgroups and
embeddable quantum homogeneous spaces is bijective only provided some
faithful flatness conditions on the module and comodule structures are 
satisfied (see \cite{MS} for more details).

Let $g:\A\rightarrow {\bf C}$ be a character, {\it i.e.} a
$*$-homomorphism of $\A$ in {\bf C}, let then 
${\rm Ad}_g\,x=\sum_{(x)}g(S(x_{(1)})) g(x_{(3)})\,x_{(2)}$. The
following Proposition is a straightforward consequence of the fact
that ${\rm Ad}_g$ is an algebra and coalgebra isomorphism and
commutes with $\tau$. 
\medskip
\begin{proposition}
\label{ad}
Let $r:\A\rightarrow\K$ be the projection that defines the right
coisotropic subgroup $\K$. Let also $r_g[x]=r[{\rm Ad}_g^{-1}\,x]$,
$x\in \A$. Then ${\rm Ker}\,r_g={\rm Ad}_g\, {\rm Ker}\,r$ is a 
$\tau$ invariant right ideal and two sided coideal and determines 
the coisotropic subgroup $\K_g$. The corresponding homogeneous 
space $B^{r_g}={\rm Ad}_g\,B^r$ is isomorphic to $B^r$ as left 
comodule algebra. \fidi
\end{proposition}   
\medskip
We remark that the two quotient structures $\K$ and $\K_g$ are 
isomorphic as coalgebras but not as right modules. Indeed 
$r_g[x\,f]=r_g[x]\,{\rm Ad}_g\,f$.

\bigskip\bigskip
\noindent
{\bf 3.} Let $q$ be a complex number of modulus one, not a root of unity.
The function algebra on the quantum group $\slqr$ is defined as 
the unital $*$-algebra generated by four real elements $a$, $b$, $c$, 
$d$ with relations
\begin{eqnarray}
a b=q\, b a\quad & a c=q\, c a \qquad b c=c b\qquad 
b d=q\, d b &\quad c d=q\, d c\nonumber\\
\label{relfq}
&d a-q^{-1}\,b c=1= a d-q\, b c\,.&
\end{eqnarray}
The Hopf algebra structure, in matrix form, is given by
\begin{equation*}
\De\left(\begin{array}{cc} a & b\\ c & d\end{array}\right)=
\left(\begin{array}{cc} a & b\\ c & d\end{array}\right)\tens
\left(\begin{array}{cc} a & b\\ c & d\end{array}\right)
\end{equation*}
\begin{equation}
\label{hopffq}
\vep \left(\begin{array}{cc} a & b\\ c & d\end{array}\right)=
\left(\begin{array}{cc} 1 & 0\\ 0 & 1\end{array}\right)\qquad
S\left(\begin{array}{cc} a & b\\ c & d\end{array}\right) =
\left(\begin{array}{cc} d & -q^{-1}\, b\\ -q\,c & a\end{array}\right)
\end{equation}
A direct application of diamond lemma leads to
\medskip
\begin{proposition}
\label{basis}
The elements 
\begin{equation}
\{ a^r b^s c^t, b^s c^t d^r\quad ; \quad r, s, t\in \Na \}
\end{equation}
form a basis of $\slqr$ as a complex vector space. \fidi
\end{proposition}

It is straightforward to verify that all the characters of $\slqr$
are of the form
\begin{eqnarray} \label{char} 
g_\alpha \, \left(\begin{array}{cc} a & b\\ c & d\end{array}\right) &
= & \left(\begin{array}{cc} \alpha & 0\\ 0 &1/\alpha\end{array}\right),
\end{eqnarray} 
with $\alpha \in {\bf R}\backslash \{0\}$.

\bigskip

We now describe the family of quantum coisotropic subgroups of $\slqr$.
 
\medskip
\begin{proposition}
Let $\C_{\mu\nu}$ be the linear subspace of $\slqr$ spanned by 
$\{(a-d + 2q^\half\,\mu\,b)\,,(\bcn)\}$ with $\mu\,,\; \nu \in {\bf R}$. 
Then $\C_{\mu\nu}$ is a $\tau$-invariant two sided coideal in $\slqr$.
\end{proposition}

{\it Proof}.
From the homomorphism property $\tau(ab)=\tau(a)\,\tau(b)$, the
$\tau$-invariance of $\C_{\mu\nu}$ is proved verifying that
\begin{equation*}
\tau(a-d + 2q^\half\,\mu\,b)=-(a-d + 2q^\half\,\mu\,b)\quad 
\mbox{and}\quad\tau(\bcn)=-q^{-1}\,(\bcn)\;.
\end{equation*}
From
\begin{equation*}
\vep(a-d)=\vep(b)=\vep(c)=0\;,
\end{equation*}
\begin{align*}
\De(a-d+2q^\half\,\mu\,b)&=(a-d+2q^\half\,\mu\,b)\tens (a+2 
q^\half\,\mu\,b)+ (d-2q^\half\,\mu\,b)\tens\\&\quad\; (a-d
+2q^\half\,\mu\,b)+ b\tens (\bcn)-(\bcn)\tens b
\intertext{and}
\De(\bcn)&=(a+2q^\half\,\mu\,b)\tens (\bcn)+(\bcn)\tens
(d-2q^\half\,\mu\,b)\\
&\quad-(a-d+2q^\half\,\mu\,b)\tens c+c\tens (a-d+2q^\half\,\mu\,b)
\end{align*}
it follows that $\C_{\mu\nu}$ is a two-sided coideal. \fidi

\medskip
Let now $\R_{\mu\nu} = \C_{\mu\nu}\; \slqr$ and $\L_{\mu\nu} = 
\slqr\;\C_{\mu\nu}$ the right and left ideals generated by
$\C_{\mu\nu}$; we call $\K_{\mu\nu}=\slqr/\R_{\mu\nu}\,$,
${}_{\mu\nu}\K=\slqr/\L_{\mu\nu}$ and
$r_{\mu\nu}$ and $\ell_{\mu\nu}$ the corresponding quotient morphisms.
We denote by $\cdot$ the action of $\slqr$ on the quotients.

In the following of this section we will study the coisotropic 
quantum subgroups defined by the right projection.

\medskip
\begin{lemma}
\label{subgroup}
Let $w_0=1$ and  
\begin{equation*}
w_n=(a+q^\half\,\chi_\sigma\,b)\dotsm (a+q^{|n|-\half}\,
\chi_\sigma\, b)\,,\quad n\in\Ze\backslash \{0\}
\end{equation*}
where $\sigma=n/|n|$ and 
$$\chi_\sigma= \mu+\sigma\sqrt{\mu^2-\nu}
:=\mu+\sigma\,\exp{\{i\pi(1-{\rm sign}(\mu^2-\nu))/4}\}\,\sqrt{|\mu^2-\nu|}\,.$$
Then $v_n = r_{\mu\nu}[w_n]$ is group-like, $\tau(v_n)=v_n\;$ 
for $\mu^2<\nu$ while $\tau(v_n)=v_{-n}$ for $\mu^2\geq\nu$. 
\end{lemma}

{\it Proof}.
From the relations $w_n\,(a+q^{-n+\half}\chi_-\,b)=(a^2+\nu\,b^2+
2q^{-\half}\mu\,ab)\,w_{n-1}$ for $n>0$ and
$w_n\,(a+q^{n+\half}\chi_+\,b)=(a^2+\nu\,b^2+
2q^{-\half}\mu\,ab)\,w_{n+1}$ 
for $n<0$ and the equality
$r_{\mu\nu}[a^2+\nu\,b^2+2q^{-\half}\mu\,ab]=r_{\mu\nu}[1]$ we have 
that
\begin{equation}
\label{recurs}
v_{n+1}=v_n\cdot (a+q^{n+\half}\chi_+\,b)\quad\mbox{and}\quad
v_{n-1}=v_n\cdot (a+q^{-n+\half}\chi_-\,b)\,, \quad n\in\Ze.
\end{equation}
We first prove, by induction, that the elements $v_n$ are group-like. 
This is trivially true for $n=0$. 
Assume then $\De v_n=v_n\tens v_n$ for $n>0$. After some algebraic 
rearrangement, making use of the induction hypothesis and the property 
\begin{equation}
\label{right_module_c}
r_{\mu\nu}[x]\cdot c=-q\,\nu\, r_{\mu\nu}[x]\cdot b \quad x\in\slqr\,,
\end{equation}
the equality $\De v_{n+1}=v_{n+1}\tens v_{n+1}$ is reduced to the
following relation 
\begin{equation}
\label{right_module_d}
v_n\cdot d=v_n\cdot(a+(q^{n+\half}\,\chi_+ 
+q^{-n+\half}\,\chi_-)\,b)\,, \quad n\in\Ze\,, 
\end{equation}
that can be proved again by induction. The results about $\tau$ are
obtained using the same procedure. The proof is similar for 
$n<0$. \fidi

\medskip
\begin{remark}
{\rm 
Analogously the left quantum subgroup obtained as the image of
$\ell_{\mu\nu}$ contains the group-like elements $\ell_{\mu\nu}[\tilde
w_n]$ where $\tilde w_0 = 1$ and 
\begin{equation}
\tilde w_n=(a+q^{|n|-\half}\,\chi_\sigma\, b)\dotsm 
(a+q^\half\,\chi_\sigma\,b)\,,\quad n\in \Ze\backslash\{0\}.
\end{equation}
}
\end{remark}

According to the structures contained in coisotropic
subgroups different notions of equivalence can be defined. In the 
following Proposition we classify them as morphisms of coalgebras 
and of modules and coalgebras.
\medskip
\begin{proposition}
\label{iso}
\begin{itemize}
\item[i{\rm \,)}] The subgroup $\K_{\mu\nu}$ is isomorphic as a
coalgebra to ${\bf R}_+$ if $\mu^2>\nu$. It is isomorphic to ${\bf S}^1$ if 
$\nu>0$ and $\mu/\sqrt{\nu}=\cos \phi_{\mu\nu}<1$ for $\cos^2\phi_{\mu\nu}\neq 
\cos^2\ell \phi$ where $\ell\in \Ze$ and 
$q=e^{\displaystyle i \phi}$. 
The coalgebras corresponding to the special series 
$\cos^2 \phi_{\mu\nu}=\cos^2 \ell \phi$ are mutually isomorphic.
\item[ii{\rm\,)}] Let $\nu>0$, $\mu/\sqrt{\nu}=\cos \phi_{\mu\nu}\leq 1\,$ and 
$\mu_n/\sqrt{\nu}=\cos(\phi_{\mu\nu}+n\, \phi)$. Then $\K_{\mu\nu}$ and 
$\K_{\mu_n\nu}$ are isomorphic as coalgebras and modules.
\end{itemize}
\end{proposition}

{\it Proof}.
We first prove the part {\it ii}\,) of the Proposition. From 
(\ref{right_module_c}), (\ref{right_module_d}) and the definition
of $C_{\mu_n\,\nu}$ it is straightforward to derive that
${\rm Ann}_{\scriptscriptstyle \K_{\mu\nu}}(v_n)={\rm
Ker}\,(r_{\mu_n\nu})$. From (\ref{recurs}) we see that 
$v_0\in v_n\cdot\slqr$ for each $n\in \Ze$. Since $\K_{\mu\nu}$ is
generated as a module by $v_0$, we conclude that
$$\K_{\mu\nu}=v_n\cdot \slqr\,\simeq\,\slqr/{\rm
Ann}_{{}_{\K_{\mu\nu}}}(v_n)
=\slqr/ {\rm Ker}\,(r_{\mu_n\nu})=\K_{\mu_n\nu}\,.$$ 
The module morphism $i:\K_{\mu\nu}\rightarrow \K_{\mu_n\nu}$ is defined by 
$i\,(v_n\cdot a)=r_{\mu_n\nu}[a]$ and it is clearly a coalgebra morphism.

We now prove the statement {\it i}\,).   
Using Proposition \ref{basis} it is easy to show that
$\K_{\mu\nu}$ is spanned by elements $r_{\mu\nu}[b^s]$ and 
$r_{\mu\nu}[ab^s]\;, s\in \Na$. From (\ref{recurs}) we get
\begin{eqnarray}
\label{vbva}
q^\half(q^n\,\chi_+ -q^{-n}\,\chi_-)\,v_n\cdot
b&=& v_{n+1}-v_{n-1}\,,
\\
(q^n\,\chi_+ - q^{-n}\,\chi_-)\,v_n\cdot a &=& q^{n}\,\chi_+
\,v_{n-1}-q^{-n}\,\chi_-\,v_{n+1}\,.\nonumber
\end{eqnarray}
Assuming  $\cos^2 \phi_{\mu\nu}\neq \cos^2 \ell \phi$, namely  
$(q^\ell\,\chi_+ - q^{-\ell}\,\chi_-)\neq 0$ for 
$\ell\in\Ze$, using a recurrence procedure starting from $s=0$ we find
\begin{equation}
r_{\mu\nu}[b^s]=\sum^s_{k=0}C^s_k\;v_{s-2k}
\end{equation}
where
\begin{equation*}
C^s_k= (-)^k\;
q^{-s/2}\;\fraz{[s]_q!}{[k]_q!\,[s-k]_q!}
\prod^s_{\substack{i=0\\i\neq s-k}}
\fraz{1}{q^{i-k}\,\chi_+ - q^{-i+k}\,\chi_-}\ \cdot
\end{equation*}
Here we use the standard notation for the $q$-numbers and $q$-factorials.
This proves, together with $r_{\mu\nu}[ab^s]=q^s\,r_{\mu\nu}[b^s]\cdot
a$ that $\K_{\mu\nu}$ is spanned by the $v_n$, with $n\in\Ze$. 
From Lemma \ref{subgroup} we have that for $\mu^2>\nu$ the subgroup
$\K_{\mu\nu}$ is isomorphic as real coalgebra to ${\bf R}_+$ and for
$\mu^2<\nu$ to ${\bf S}^1$. 

For each character $g_\alpha$ defined in {\rm (\ref{char})}, according 
to Proposition \ref{ad}, the real coisotro\-pic quantum subgroup 
$(\K_{\mu\nu})_{g_\alpha}$ is generated by 
$\C_{\mu_\alpha\nu_\alpha}={\rm Ad}_{g_\alpha}\,\C_{\mu\nu}$ 
with $\mu_\alpha=\mu/\alpha^2$ and $\nu_\alpha=\nu/\alpha^4$.
By observing that two subgroups of the special series can be connected
by the composition of the adjoint map and a morphism introduced in 
{\it ii}\,) we conclude that they are isomorphic as coalgebras.
\fidi

In the following we give an explicit description of the coalgebra
corresponding to the special series 
$\cos^2 \phi_{\mu\nu}=\cos^2 \ell \phi$. By using the result of
Proposition \ref{iso} {\it ii}\,) we can assume $\ell=0$, {\it i.e.}
$\mu^2=\nu$. We denote such coalgebra as $\Re_q$. 

\medskip
\begin{lemma}
\label{xn}
For each $n\geq 1$ let 
$$X_n=\sum_{i=1}^n \fraz{(q-q^{-1})^{i-1}\,q^{i/2}\,\mu^i\,[n-1]_q!}
      {[i]_q\,[n-i]_q!}\, v_{n-i}\cdot b^i\;.
$$
Then
$$\De X_n=X_n\tens v_n+v_n\tens X_n\;,\quad\quad \tau(X_n)=-X_n\;,$$
and $\{v_n\}_{n\in \Na}\cup \{X_n\}_{n\geq 1}$ is a linear basis for $\Re_q$.
\end{lemma}  

{\it Proof}.
The coproduct and the $\tau$ of $X_n$ can be calculated by induction 
observing that $X_1=q^\half\mu\,v_0\cdot b$ and using the relation 
$X_n\cdot (q^{-n}\,a - q^n\,d + 2 q^\half \mu\,b)=-[n+1]_q(q-q^{-1})\,X_{n+1}$.

Let's define, for each $k\geq 0$, a linear mapping 
$g_k:\;\Re_q\rightarrow \Co$ such that $g_k(v_n)=\delta_{kn}$. 
If we suppose that
$\sum_k(\alpha_k\, X_k+\beta_k\, v_k)=0$, for some
$\alpha_k\ ,\beta_k$ and we apply $({\rm id}\tens g_k)\De$ we obtain 
that $\alpha_k\,X_k\in {\rm Span}\{v_l\}_{l\in \Na}$. Since $X_n$
cannot be generated by group-like elements, we 
conclude that $\alpha_k=0$ and thus $\beta_k=0$. Therefore
$\{v_n,X_n\}$ are linearly independent.
Finally we observe that \(\Re_q={\rm Span}\{v_0\cdot b^k,v_1\cdot 
b^k\}_{k\in \Na}\). As 
\(v_1\cdot b^s\in {\rm Span}\{v_k\}_{k\leq s+1} \oplus 
{\rm Span}\{v_0\cdot b^k\}_{k\leq s-1}\)
it is sufficient to show that 
$v_0\cdot b^s\in {\rm Span}\{v_k,X_k\}_{k\leq s}$. This can be done by 
induction, when the relation 
$$q^\half \mu(q^n-q^{-n})X_n\cdot b=\fraz{[n+1]_q}{[n]_q}\,X_{n+1}-
\fraz{[n-1]_q}{[n]_q}\,X_{n-1}- q^\half \mu \fraz{(q^n+q^{-n})}{[n]_q}
\,v_n\cdot b\,$$
is used. \fidi

In the following proposition we summarize the properties of the
coalgebra $\Re_q$.
\medskip
\begin{proposition}
\label{Re_q}
\begin{itemize}
\item[i{\rm \,)}] If $\Re_q^{(n)}={\rm Span}\{v_n,X_n\}$ for any fixed 
$n>0$ and $\Re_q^{(0)}={\rm Span}\{v_0\}$, then
$\Re_q=\bigoplus_{n\in\Na}\,\Re_q^{(n)}$ as a direct sum of
cocommutative coalgebras.
\item[ii{\rm\,)}] The only simple subcoalgebras are those generated 
by $v_n$, i.e. $\Re_q$ is pointed. The corresponding one dimensional 
corepresentations are unitary.
\item[iii{\rm\,)}] There exists no Hopf algebra isomorphic to
$\Re_q$ as coalgebra. 
\end{itemize}
\end{proposition}

{\it Proof}.
The point {\it i\,}) is a direct consequence of Lemma
(\ref{xn}).

Using Theorem (8.0.3) of \cite{Swee} we have that each simple
coalgebra must be contained in $\Re_q^{(n)}$ for some $n$. It is clear
that the only simple subcoalgebra of $\Re_q^{(n)}$ is ${\rm
Span}\{v_n\}$. This proves the point {\it ii\,}). 

Finally let us suppose that $\Re_q$ has a bialgebra structure and 
let $v_{n_0}=1$. By Theorem (8.1.1) of \cite{Swee} the irreducible 
component $\Re_q^{(n_0)}$ of $v_{n_0}$ must be a sub-bialgebra. 
It is not difficult to see that the only two dimensional 
bialgebra is generated by two group-like elements, so that 
$\Re_q^{(n_0)}$ must be one dimensional, {\it i.e.} $n_0=0$.
By using Theorem (8.1.5) of \cite{Swee} if the irreducible component
of the identity is one dimensional then $\Re_q$ must be isomorphic as
an algebra to ${\rm Span}\{v_n\}_{n\in\Na}$. This is not true and
shows point {\it iii\,}). \fidi

\medskip
\begin{remark}
{\rm
{\it Classical limit of the special series.} 
The equation (\ref{vbva}) for $\mu^2=\nu$ reads
$\;v_{n+1}-v_{n-1}=(q-q^{-1})\,q^\half\,\mu\; [n]_q\;v_n\cdot
b$. As a consequence we have that
$\;\{v_{2n}-v_0\,,\,v_{2n+1}-v_1\}\in(q-q^{-1})\,\Re_q\,$, so that 
in the classical limit all the group-like elements $v_n$ collapse 
into $v_0$ or $v_1$ according to the parity of $n$.
}
\end{remark}

\bigskip\bigskip
\noindent
{\bf 4.} Starting from a real coisotropic quantum subgroup one has a 
naturally defined real embeddable quantum homogeneous space
\begin{equation*}
B^{r_{\mu\nu}}=\{ x\in\slqr\, : \, (id\tens r_{\mu\nu})\De x=x\tens 
r_{\mu\nu}[1]\}\,.
\end{equation*}

Let 
\begin{align}
\label{z1z2z3}
z_1&=q^{-\half}(ac+\nu\,bd)+2 \mu\,bc\,,\notag \\
z_2&=c^2+\nu\, d^2+2\mu q^{-\half} cd\,,\\
z_3&=a^2+\nu\, b^2+2\mu q^{-\half} ab\,.\notag
\end{align}
By direct computation it can be verified that $z_i$ are real elements 
with relations
\begin{align}
\label{relz1z2z3}
z_1z_2&=q^2 z_2z_1\,,
\notag \\
z_1z_3&=q^{-2}z_3z_1\,,\\
z_3z_2&=\nu+q^2z_1^2+2 \mu q\,z_1 \notag
\end{align}
and  coproducts
\begin{equation}
\begin{split}
\De(z_1)= &\,(1+(q+q^{-1})bc)\tens z_1 + q^{-\half}bd\tens z_2 + 
q^{-\half}ac\tens z_3\\ 
&+2 \mu\,bc\tens 1\,,\\
\De(z_2)= &\,q^{-\half}(q+q^{-1})\,cd\tens z_1 + d^2\tens z_2 +
c^2\tens z_3 \\
&+ 2 \mu q^{-\half}\,cd\tens 1\,,\\
\De(z_3)= &\,q^{-\half}(q+q^{-1})\,ab\tens z_1 + b^2\tens z_2 + a^2\tens
z_3\\
&+2 \mu q^{-\half}\,ab\tens 1\,.
\end{split}
\end{equation}

The following proposition can be proven according to the lines
suggested in Proposition 5.4 of \cite{Br}.

\medskip
\begin{proposition}
The left comodule subalgebra $B^{r_{\mu\nu}}$ is generated by 
$\{z_i\}_{i=1,2,3}$. \fidi
\end{proposition}
\medskip
The characters of $B^{r_{\mu\nu}}$ are given by
$g(z_1)=0\,,\,g(z_2)=\alpha \nu\,,\,g(z_3)=1/\alpha$ with 
$\alpha\in\Re\backslash\{0\}$. Using the map 
$({\rm id}\otimes g)\circ \Delta$ we see that $B^{r_{\mu\nu}}$ 
and $B^{r_{\mu'\nu'}}$ are isomorphic as left comodule algebras 
if $\mu'=\mu/\alpha$ and $\nu'=\nu/\alpha^2$ with 
$\alpha\in \Re\backslash\{0\}$. We remark that $g$ is a restriction 
of a character defined on the whole $\slqr$ only for $\alpha>0$.

\begin{remark}
{\rm
In the classical limit the last of the relations (\ref{relz1z2z3}) 
reads $z_3z_2=\nu+z_1^2+2\mu z_1$; posing $z_1=z-\mu$, $z_2=x+y$, 
$z_3=x-y$ we get $x^2-y^2-z^2=\nu-\mu^2$. Therefore $\Theta=\nu-\mu^2$ 
is the parameter, defined in the Introduction, that classifies the
homogeneous spaces of $\slr$. 
}
\end{remark}

For all $\mu\,,\nu$ the irreducible corepresentations of the subgroup are one
dimensional and defined by  $\rho_j(1)=v_j$ where $j\in\Ze$ for 
$\Re_+\,,{\bf S}^1$ and $j\in\Na$ for $\Re_q$. 
We remark that they are unitary only in the case of ${\bf S}^1$ and $\Re_q$.

Following the scheme of \cite{BCGST} and \cite{Ci2}, we induce
respectively right and left corepresentations of the whole 
quantum group $\slqr$ on
\begin{eqnarray*}
B^{r_{\mu\nu}}_j=\{ x\in\slqr \big| (id\tens r_{\mu\nu})\De x=x\tens v_j\}\,\\
{}_jB^{\ell_{\mu\nu}}=\{x\in\slqr \big| (\ell_{\mu\nu}\tens id)\De x=v_j\tens
x\}\;.
\end{eqnarray*}
The coaction map is simply the restriction of the coproduct to
these spaces.

As \((\ell_{\mu\nu}\tens {\rm id})\De:
B^{r_{\mu\nu}}_j\rightarrow \K_{\mu\nu}\tens B^{r_{\mu\nu}}_j\) we can
also define the left and right corepresentations 
${}_jB^{\mu\nu}_k={}_kB^{\ell_{\mu\nu}}\cap B^{r_{\mu\nu}}_j$. 
By direct computation we can give the explicit characterization of the 
double coset
$${}_oB^{\mu\nu}_o ={\rm Span}\{(z_2+\nu\,z_3-2\mu\,z_1)^n\}_{n\in\Na}.$$

Since ${\bf S}^1$ and $\Re_+$ are cosemisimple we can apply Corollary
1.5 from \cite{MS}; we then have, for the corresponding values of $\mu$ 
and $\nu$, the following decomposition
$$\slqr=\bigoplus_{j\in\Ze} {}_jB^{\ell_{\mu\nu}} =\bigoplus_{j\in\Ze} 
S(B^{r_{\mu\nu}}_j) \;.$$
Furthermore $B^{r_{\mu\nu}}_j$ (${}_jB^{\ell_{\mu\nu}}$) is a finitely 
generated projective
$B^{r_{\mu\nu}}$ ($B_{\ell_{\mu\nu}}$)--module. This property is
usually rephrased by saying that $B^{r_{\mu\nu}}_j$
(${}_jB^{\ell_{\mu\nu}}$) 
is the space of sections of a quantum line bundle on the quantum 
homogeneous space. 
The space further decomposes as 
$$\slqr = \bigoplus_{jk}{}_jB^{\mu\nu}_k \,.$$ 
The study of these spaces will be given in a forthcoming paper
(see \cite{BrMj} for $B_1$ for the case of the quantum spheres). 
We remark, however, that no conclusion can be drawn on vector bundles in the
case of $\Re_q$ {\it via} ref. \cite{MS}.

\end{document}